\documentclass[reqno,draft]{amsart}

\usepackage{amsthm}
\usepackage{amssymb}
\usepackage{amsmath}

\newcommand{\labell}[1] {\label{#1}}

\numberwithin{equation}{section}

\newtheorem {Theorem}   {Theorem} 
\numberwithin{Theorem}{section}

\newtheorem {Lemma}[Theorem]    {Lemma}         
\newtheorem {Proposition}[Theorem]{Proposition}  
\theoremstyle{definition}

\theoremstyle{remark}
\newtheorem{Remark}[Theorem]{Remark}

\newcounter{cond}

\expandafter\chardef\csname pre amssym.def at\endcsname=\the\catcode`\@ 
\catcode`\@=11 
\def\undefine#1{\let#1\undefined} 
\def\newsymbol#1#2#3#4#5{\let\next@\relax 
 \ifnum#2=\@ne\let\next@\msafam@\else 
 \ifnum#2=\tw@\let\next@\msbfam@\fi\fi 
 \mathchardef#1="#3\next@#4#5}
\def\mathhexbox@#1#2#3{\relax 
 \ifmmode\mathpalette{}{\m@th\mathchar"#1#2#3}%
 \else\leavevmode\hbox{$\m@th\mathchar"#1#2#3$}\fi} 
\def\hexnumber@#1{\ifcase#1 0\or 1\or 2\or 3\or 4\or 5\or 6\or 7\or 8\or 
 9\or A\or B\or C\or D\or E\or F\fi} 
 
\font\teneufm=eufm10 
\font\seveneufm=eufm7
\font\fiveeufm=eufm5 
\newfam\eufmfam
\textfont\eufmfam=\teneufm 
\scriptfont\eufmfam=\seveneufm 
\scriptscriptfont\eufmfam=\fiveeufm 
 
\catcode`\@=\csname pre amssym.def at\endcsname 

\def    \eps    {\epsilon}
\def    \p      {\partial}
\def    \px     {\partial_x}
\def    \py     {\partial_y}
\def    \pt     {\partial_t}
\def    \pxi    {\partial_\xi}
\def    \peta   {\partial_\eta}
\def    \pu     {\partial_u}

\def    \lt     {\left(}
\def    \rt     {\right)}

\newcommand{\CI}{{\mathcal J}}

\newcommand{\const}{{\mathit const}}

\newcommand{\D}{{\mathfrak D}}

\newcommand{\tF}{\tilde{F}}
\newcommand{\fr}{\frac}
\newcommand{\lm}{\displaystyle\lim}

\def    \reals  {{\mathbb R}}
\def    \R      {{\mathbb R}}
\def    \integers       {{\mathbb Z}}
\def    \Z      {{\mathbb Z}}

\def    \p      {\partial}

\def    \ssminus        {\smallsetminus}

\def    \inv    {{\operatorname{inv}}}
\def    \Pin    {\Phi_y^{\inv}}

\begin{document}


\setlength{\smallskipamount}{6pt}
\setlength{\medskipamount}{10pt}
\setlength{\bigskipamount}{16pt}





\title[A $C^2$-counterexample to the Hamiltonian Seifert Conjecture]{A 
$C^2$-smooth counterexample to the Hamiltonian Seifert Conjecture in $\R^4$}

\author[Viktor Ginzburg]{Viktor L. Ginzburg}
\author[Ba\c sak G\"urel]{Ba\c sak Z. G\"urel}

\address{Department of Mathematics, UC Santa Cruz, 
Santa Cruz, CA 95064, USA}
\email{ginzburg@math.ucsc.edu, basak@math.ucsc.edu}

\date{\today}

\thanks{The work is partially supported by the NSF and by the faculty
research funds of the University of California, Santa Cruz.}

\bigskip

\begin{abstract}
We construct a proper $C^2$-smooth function on $\R^4$ such that its 
Hamiltonian flow has no periodic orbits on at least one regular level set. 
This result 
can be viewed as a $C^2$-smooth counterexample to the Hamiltonian Seifert 
conjecture in dimension four. 
\end{abstract}

\maketitle

\section{Introduction} 
The ``Hamiltonian Seifert conjecture'' is the question whether or not there
exists a proper function on $\R^{2n}$ whose Hamiltonian flow has no periodic
orbits on at least one regular level set. We construct a $C^2$-smooth function
on $\R^4$ with such a level set. Following the tradition 
of \cite{kug,kugk,kuk,kuk:icm,schweitzer}, this result 
can be called a $C^2$-smooth counterexample to the Hamiltonian Seifert 
conjecture in dimension four. 

In dimensions greater than six, $C^\infty$-smooth counterexamples to 
the Hamiltonian Seifert conjecture were constructed by one of the authors, 
\cite{gi:seifert}, and simultaneously by M. Herman, 
\cite{herman-fax,herman}. In dimension six, a $C^{2+\eps}$-smooth 
counterexample was found by M. Herman, \cite{herman-fax,herman}. This
smoothness constraint was later relaxed to $C^\infty$ in
\cite{gi:seifert97}. A very simple and elegant construction of a new
$C^\infty$-smooth counterexample in dimensions greater than four was
recently discovered by Kerman, \cite{ely:example}. The flow in
Kerman's example has dynamics different from the 
ones in \cite{gi:seifert,gi:seifert97,herman-fax,herman}.
We refer the reader
to \cite{gi:bayarea,gi:barcelona} for a detailed discussion of the 
Hamiltonian Seifert conjecture. The reader interested in the results
concerning the original Seifert conjecture settled by K. Kuperberg,
\cite{kugk,kuk}, should consult \cite{kuk:icm,kuk:notices}. Here we only 
mention that a $C^1$-smooth counterexample to the Seifert conjecture
on $S^3$ was constructed by Schweitzer, \cite{schweitzer}, and 
a $C^1$-smooth volume--preserving counterexample  
on $S^3$ was found by G. Kuperberg, \cite{kug}. The ideas from both
of these constructions play an important role in this paper.

An essential difference of the Hamiltonian case from the general one
is manifested by the almost existence theorem, 
\cite{ho-ze:per-sol,ho-ze:book,str}, which asserts that almost all regular
levels of a proper Hamiltonian have periodic orbits (see Remark 
\ref{rmk:almost-ex}). In other words, regular levels 
without periodic orbits are exceptional in the sense of measure theory.

The existence of a $C^2$-counterexample to the Hamiltonian Seifert conjecture
in dimension four was announced by the authors in \cite{GG}, where a
proof was also outlined.
Here we give a detailed construction of this counterexample.

\subsection*{Acknowledgments.} The authors are deeply grateful to 
Helmut Hofer, Anatole Katok, Ely Kerman, Krystyna Kuperberg, 
Mark Levi, Debra Lewis, Rafael de la Llave, Eric Matsui, and Maria Schonbek 
for useful discussions and suggestions.

\section{Main Results}
\labell{sec:main}

Recall that \emph{characteristics} on a hypersurface $M$ in a
symplectic manifold $(W,\eta)$ are, by definition, the (unparameterized) 
integral curves of the field of directions $\ker (\eta|_M)$. 

Let $\R^{2n}$ be equipped with its standard symplectic structure.

\begin{Theorem}
\labell{thm:main}
There exists a $C^2$-smooth embedding $S^3\hookrightarrow\R^4$ which
has no closed characteristics. This embedding
can be chosen $C^0$-close and $C^2$-isotopic to an ellipsoid.
\end{Theorem}

As an immediate consequence we obtain

\begin{Theorem}
\labell{thm:main2}
There exists a proper $C^2$-function $F\colon\R^4\to \R$ such that
the level $\{F=1\}$ is regular and the Hamiltonian flow of $F$
has no periodic orbits on $\{F=1\}$. In addition, 
$F$ can be chosen so that this level is $C^0$-close and $C^2$-isotopic to an
ellipsoid. 
\end{Theorem}

\begin{Remark}
\labell{rmk:almost-ex}
Regular levels of $F$ without periodic orbits are exceptional in the
sense that the set of corresponding values of $F$ has zero measure. This is a 
consequence
of the almost existence theorem, \cite{ho-ze:per-sol,ho-ze:book,str},
which guarantees that for a $C^2$-smooth (and probably even $C^1$-smooth)
function, periodic orbits exist
on a full measure subset of the set of regular values. In particular,
since all values of $F$ near $F=1$ are regular, almost all levels of
$F$ near this level carry periodic orbits.
\end{Remark}

\begin{Remark}
It is quite likely that our construction gives an 
embedding $S^3\hookrightarrow\R^4$ without closed characteristics,
which is \emph{$C^{2+\alpha}$-smooth}.
\end{Remark}

\begin{Remark}
\labell{rmk:main}
Similarly to its higher-dimensional counterparts, 
\cite{gi:seifert,gi:seifert97},
Theorem \ref{thm:main} extends to other symplectic manifolds as follows.
Let $(W,\eta)$ be a four-dimensional symplectic manifold and let
$i\colon M\hookrightarrow W$ be a $C^\infty$-smooth embedding such that 
$i^*\eta$ has only a 
finite number of closed characteristics. Then there exists a $C^2$-smooth 
embedding $i'\colon M\hookrightarrow W$, which is $C^0$-close and isotopic 
to $i$, such that ${i'}^*\eta$ has no closed characteristics.
\end{Remark}

The rest of the paper is devoted to the proof of Theorem \ref{thm:main}. 
The idea of the proof is to adjust 
Schweitzer's construction, \cite{schweitzer}, of an aperiodic $C^1$-flow 
on $S^3$ to make it embeddable into $\R^4$ as a Hamiltonian
flow. This is done by introducing a Hamiltonian version of Schweitzer's 
plug. More specifically, the flow on Schweitzer's plug is defined as the 
Hamiltonian flow of a certain multi-valued function $K$ 
which we use to find a symplectic embedding of the plug
(see Proposition \ref{prop:K1} and Remark \ref{rmk:j}). The existence of such
a function $K$ heavily depends on the choice of a Denjoy vector field in
Schweitzer's plug. Namely, the Denjoy vector field is required
to be essentially as smooth as a Denjoy vector field can be (see Remark 
\ref{rmk:smooth}). Implicitly,
the idea to define the flow on Schweitzer's plug using the Hamilton equation
goes back to G. Kuperberg's paper \cite{kug}. 

As of this moment we do not know if G. Kuperberg's flow can be embedded into
$\R^4$. The two constructions differ in an essential way. 
The Denjoy flow  and the function $K$ in G. Kuperberg's example are 
required to have properties very different from the ones we need.
As a consequence, our method to embed the plug into
$\R^4$ does not apply to G. Kuperberg's plug.
(For example, one technical but essential discrepancy
between the methods is as follows. In G. Kuperberg's 
construction, it is important to take a rotation number which cannot be
too rapidly approximated by rationals, while the Denjoy map is not required to
be smoother than just $C^1$. On the other hand, in our construction the value
of rotation number is irrelevant, but the smoothness of the Denjoy map
plays a crucial role.)

The proof is organized as follows. In Section 
\ref{sec:embedding} we describe the symplectic embedding of
Schweitzer's flow assuming the existence of 
the plug with required properties. In Sections \ref{sec:pfK1} and
\ref{sec:proofs} we derive the existence of such a flow on the plug
from the fact (Lemma \ref{lemma:main}) that there exists a 
``sufficiently smooth'' Denjoy flow
on $T^2$. Finally, this ``sufficiently smooth'' Denjoy flow is 
constructed in Section \ref{sec:pflemma-main}.

\section{Proof of Theorem \ref{thm:main}: The Symplectic Embedding}
\labell{sec:embedding}

Let us first fix some notations. Throughout this paper $\sigma$ 
denotes the standard symplectic form on $\R^{2m}$ or the pull-back of
this form to $\R^{2m+1}$ by the projection $\R^{2m+1}\to \R^{2m}$ along
the first coordinate; $I^{2m}$ stands for a cube in $\R^{2m}$ whose
edges are parallel to the coordinate axis. The product $[a,b]\times I^{2m}$ is
always assumed to be embedded into $\R^{2m+1}$ (henceforth, the standard 
embedding) so that the interval $[a,b]$ is 
parallel to the first coordinate. We refer to the direction along the
first coordinate $t$ (time) in $\R^{2m+1}$ (or $[a,b]$ in 
$[a,b]\times I^{2m}$) as the vertical direction.

All maps whose smoothness is not specified are $C^\infty$-smooth.

Theorem \ref{thm:main}, as do similar theorems in dimensions greater
than four, follows from the existence of a symplectic plug. The definitions
of a plug vary considerably (see \cite{gi:seifert,ely:example,kug}), and here 
we use the one more suitable for our purposes.

A \emph{$C^k$-smooth symplectic plug} in dimension $2n$ is a $C^k$-embedding 
$J$ of $P=[a,b]\times I^{2n-2}$ into $P\times \R\subset\R^{2n}$ such that

\begin{enumerate}

\item[P1.] 
\emph{The boundary condition}: The embedding $J$ is the identity
embedding of $P$ into $\R^{2n-1}$ near the boundary $\partial P$. Thus
the characteristics of $J^*\sigma$ are parallel to the vertical direction
near $\partial P$.

\item[P2.]
\emph{Aperiodicity}:
The characteristic foliation of $J^*\sigma$ is aperiodic, i.e., $J^*\sigma$
has no closed characteristics.

\item[P3.] 
\emph{Existence of trapped trajectories}:
There is a characteristic of $J^*\sigma$ beginning on $\{a\}\times I^{2n-2}$ 
that never exits the plug. Such a characteristic is said to be trapped in $P$.

\item[P4.]
The embedding $J$ is $C^0$-close to the standard embedding and $C^k$-isotopic
to it.

\item[P5.]
\emph{Matched ends or the entrance--exit condition}:
If two points $(a, x)$, the ``entrance'', and $(b,y)$, 
the ``exit'', are on the same characteristic, then $x=y$. In other words,
for a characteristic that meets both the bottom and the top of the plug,
its top end lies exactly above the bottom end.

\end{enumerate}

\begin{Theorem}
\labell{thm:plug}
In dimension four, there exists a $C^2$-smooth symplectic plug.
\end{Theorem}

\begin{proof}[Proof of Theorem \ref{thm:main}]
Theorem \ref{thm:main} readily follows from Theorem \ref{thm:plug}. Consider
an irrational ellipsoid in $\R^4$ and pick two little balls each of which is
centered at a point on a closed characteristic on the ellipsoid. 
Intersections of these balls with the ellipsoid can be viewed symplectically
as open subsets in $\R^3$. By scaling the plug we can assume that 
$[a,b]\times I^2$ can be embedded into each of these open balls so that the
closed characteristic on an ellipsoid matches a trapped trajectory 
in the plug. Now we perturb the ellipsoid by means of the embedding $J$ 
within each of these open subsets. The resulting embedding has no closed
characteristics, $C^0$-close to the ellipsoid and $C^2$-isotopic to it.
\end{proof}

\begin{proof}[Proof of Theorem \ref{thm:plug}]
First observe that it suffices to construct a semi-plug, i.e., a ``plug'' 
satisfying only the conditions (P1)-(P4). Indeed, a plug can then be obtained
by combining two symmetric semi-plugs. More precisely, suppose that a 
semi-plug with embedding $J_{_-}$ has been constructed. Without loss of 
generality
we may assume that $[a,b]=[-1,0]$. Define a semi-plug on 
$[0,1]\times I^2$ with embedding $J_{+}$ by setting $J_{+}(t,x)=RJ_{-}(-t,x)$,
$t\in [0,1]$ and $x\in I^2$, where $R$ is the reflection of $\R^4$ in $\R^3$.
Combined together, these semi-plugs give rise to a plug on $[-1,1]\times I^2$.

We will construct a semi-plug by perturbing the standard embedding 
of $[a,b]\times I^2$ on a subset $M\subset [a,b]\times I^2$. This
subset is diffeomorphic to $[-1,1]\times \Sigma$, where $\Sigma$ is
a punctured torus.

It is more convenient to do this perturbation using slightly
different ``coordinates'' on a neighborhood of $M$. 
More specifically, we will first consider an embedding of $M$ into another 
four-dimensional symplectic manifold $(W, \sigma_W)$ such that
the pull-back of $\sigma_W$ is still $\sigma|_M$.
Then we $C^0$-perturb this embedding so that the characteristic vector field of
the new pull-back will have properties similar to those of Schweitzer's
plug. By the symplectic neighborhood theorem, a neighborhood of $M$ in
$W$ is symplectomorphic to that of $M$ in $\R^4$. This will allow us to
turn the embedding $M\hookrightarrow W$ into the required embedding
$J\colon M\hookrightarrow \R^4$. (See the diagrams \eqref{eq:diag1}
and \eqref{eq:diag2} below.)

To construct the perturbed embedding $M\hookrightarrow W$, we first
embed $M$ into $[-1,1]\times T^2$ by puncturing the torus in a suitable
way. Then we find a map $j\colon [-1,1]\times T^2 \to W$ such
that the characteristic vector field of $j^*\sigma_W$ is aperiodic and
has trapped trajectories. 

The embedding $j$ is constructed as follows.
Let $(x,y)$ be coordinates
on $T^2$. Consider the product $W=(-2,2)\times S^1\times T^2$ with 
coordinates $(t, x,u,y)$ and symplectic form 
$\sigma_W=dt\wedge dx+du\wedge dy$. The map $j$ is
a $C^0$-small perturbation of 
$$
j_0\colon [-1,1]\times T^2\to W;
\quad j_0(t,x,y)=(t,x,x,y).
$$
Note that $j_0(t,x,y)=(t, x,K_0,y)$, where $K_0(t,x,y)=x$. To define
$j$, let us replace $K_0$ by 
a mapping $K\colon [-1,1]\times T^2\to S^1$ to be specified later on.
In other words, set
$$
j\colon [-1,1]\times T^2\to (-2,2)\times S^1\times T^2,
\quad\text{where}\quad j(t,x,y)=(t,x,K,y).
$$
It is clear that $j$ is an embedding. (An explanation of the origin 
of $j$ is given in Remark \ref{rmk:j}.) 
The pull-back $j^*\sigma_W$ is the form 
$$
j^*\sigma_W=dt\wedge dx+ (\px K)dx\wedge dy + (\pt K)dt\wedge dy 
$$
with characteristic vector field
$$
v=(\px K)\pt-(\pt K) \px + \py.
$$
To ensure that
(P1)-(P4) hold we need to impose some requirements on $K$.

To specify these requirements, consider a Denjoy vector field 
$\py+h\px$ on $T^2$.
This vector field should satisfy certain additional conditions which will be
detailed in Section \ref{sec:pflemma-main}. 
Denote by $\D$ the Denjoy continuum for this field.\footnote{We refer the 
reader to \cite{HS,KH,schweitzer} for a discussion of
Denjoy maps and vector fields.} 
Pick a point $(x_0,y_0)$ in the complement of $\D$.
Fix a small, disjoint from $\D$, neighborhood $V_0$ of $(x_0,y_0)$. 
Consider the tubular neighborhood of the line $(t,x_0,y_0+t)$
in $[-1,1]\times T^2$ of the form 
$\{ (t,x,y+t)\mid (x,y)\in V_0,~ t\in [-1,1]\}$.
Fix also a small neighborhood of the boundary
$\p ([-1,1]\times T^2)$ and denote by $N$ the union of these 
neighborhoods.

\begin{Proposition}
\labell{prop:K1}
There exists a $C^2$-smooth mapping $K\colon [-1,1]\times T^2\to S^1$ such
that 

\begin{enumerate}

\item[K1.] $v$ is equal to the Denjoy vector field (i.e.,
$\px K =0$ and $\pt K=-h$) at every point of $\{0\}\times \D$;

\item[K2.] the $t$-component of $v$ is positive (i.e., $\px K> 0$) 
on the complement of $\{0\}\times \D$;

\item[K3.]$K$ is $C^0$-close\footnote{More specifically, 
for any $\eps>0$ there exists $K$ satisfying (K1)-(K2) and (K4) such that
$\parallel K-K_0\parallel <\eps$. The required value of $\eps$ is 
determined by the size of the neighborhood $U$ in the symplectic 
neighborhood theorem; see below.}  to the map 
$K_0\colon (t,x,y)\mapsto x$;

\item[K4.]$K=K_0$ on $N$.
\end{enumerate}
\end{Proposition}

Let us defer the proof of the proposition to Section \ref{sec:pfK1} and 
finish the proof of Theorem \ref{thm:plug}. From now on we assume that $K$ 
is as in Proposition \ref{prop:K1}.

By (K1) and (K2), $v$ has a trapped trajectory and is aperiodic.
Indeed, by (K1), $\{0\}\times \D$ is invariant under the flow of $v$ and
on this set the flow is a Denjoy flow. By (K2), the vertical 
component of $v$ is non-zero unless the point is in $\{0\}\times \D$. 
This implies that periodic orbits can only occur within $\{0\}\times \D$.
Since the Denjoy flow is aperiodic, so is the entire flow of $v$.
Furthermore, it is easy to see that since $\{0\}\times \D$ is invariant, there
must be a trapped trajectory. Furthermore, $v=\pt+\py$ on $N$ by (K4).

Now we are in a position to define $J$.
Let $\Sigma$ be the torus $T^2$, punctured at 
$(x_0,y_0)$. To be more accurate, $\Sigma$ is obtained by deleting
a neighborhood of $(x_0,y_0)$, contained in $V_0$. There exists a 
symplectic bridge 
immersion of $(\Sigma, dx\wedge dy)$ into some cube $I^2$ with the standard 
symplectic structure. Hence, there exists an embedding
$$
M=[-1,1]\times \Sigma \hookrightarrow [a,b]\times I^2\subset \R^3
\subset \R^4
$$
such that the pull back of $\sigma$ is $dx\wedge dy$. Henceforth, we 
identify $M$ with its image in $\R^4$.

On the other hand, we can embed $M$ into $[-1,1]\times T^2$ by means of
$$
\varphi\colon M=[-1,1]\times \Sigma \to [-1,1]\times T^2;
\quad \varphi(t,x,y)=(t,x,y+t).
$$
Then $\varphi_*\pt=\pt+\py$ and $(j_0\varphi)^*\sigma_W= dx\wedge dy$.
The argument similar to the proof of 
the symplectic neighborhood theorem, \cite[Lemma 3.14]{mcduff-sal},
shows (see \cite[Section 4]{gi:seifert} for details)
that a ``neighborhood'' of $M$ in $\R^4$ is
symplectomorphic to a ``neighborhood'' $U$ of $j_0\varphi(M)$ 
in $W$. More precisely, for a small $\delta>0$, there exists a 
symplectomorphism
$$
\psi\colon M\times (-\delta,\delta)\to U\subset W
$$
extending $j_0\varphi$, i.e., such that $\psi|_{M}=j_0\varphi$. These maps
form the following diagram:
\begin{equation}
\labell{eq:diag1}
\begin{array}{rrcrr}
M & \hookrightarrow &  (-\delta,\delta)\times M &  \subset & \R^4 \\
\parallel & & \downarrow \scriptstyle{\psi} & &  \\
M & {\buildrel j_0\varphi \over \longrightarrow} &  U &  \subset & W
          \end{array}
\end{equation}

By (K3),  $j$ is $C^0$-close to $j_0$. Furthermore,
$j=j_0$ on $N$ by (K4). Hence, $j$ can be assumed to take values
in $U$ (see Remark \ref{rmk:U}).

Finally, set
$$
J=\psi^{-1}j\varphi
$$
on $M$. In other words, $J$ is defined by the diagram:
\begin{equation}
\labell{eq:diag2}
\begin{array}{rrcrr}
M & {\buildrel J \over \longrightarrow}&  (-\delta,\delta)\times M &  \subset & \R^4 \\
\parallel & & \downarrow \scriptstyle{\psi} & &  \\
M & {\buildrel j\varphi \over \longrightarrow} &  U &  \subset & W
          \end{array}
\end{equation}
Then $(J^*\sigma)|_M=(j\varphi)^*\sigma_W$. To finish the definition of $J$, 
we extend it as the standard embedding to $[a,b]\times I^2\ssminus M$.

The characteristic vector field of $J^*\sigma$ is $\pt$ in the
complement of $M$ and $(\varphi^{-1})_*v$ on $M$. Since
$(\varphi^{-1})_*v=\pt$ near $\p M$, these vector fields match smoothly
at $\p M$. It is clear that (P1) is satisfied. Since $v$ has a trapped 
trajectory and is aperiodic, the same is true for $(\varphi^{-1})_*v$, i.e., 
the conditions (P2) and (P3) are met.  The condition (P4) is easy to 
verify. Hence, $J$ is indeed a semi-plug. 
\end{proof}

\begin{Remark}
\labell{rmk:U}
The following argument shows in more detail why $j$ can be assumed to take 
values in
$U$. Let us slightly shrink $M$ by enlarging the puncture in $T^2$
and shortening the interval $[-1, 1]$. Denote the resulting manifold
with corners by $M'$. The shrinking is made so that $\p M'\subset N$
and hence $M\ssminus M' \subset N$. It follows
that $U$ contains a genuine neighborhood $U'$ of $j_0\varphi(M')$. Thus, if
$K$ is sufficiently $C^0$-close to $K_0$, we have 
$j(\varphi(M'))\subset U'$. On $j_0\varphi (M\ssminus M')$, we have $K=K_0$ by
(K5) and hence $j=j_0$. Therefore, $j(\varphi(M))\subset U$.
\end{Remark}

\begin{Remark}
\labell{rmk:j}
The definition of the embedding $j$ can be explained as follows. Let us view
the annulus $[-1,1]\times S^1$ with symplectic form $dt\wedge dx$ as a 
symplectic manifold and the product $[-1,1]\times T^2$ as the extended
phase space with the $y$-coordinate being the time-variable. Then
we can regard $K$ as a (multi-valued) time-dependent Hamiltonian on 
$[-1,1]\times S^1$. The embeddings $j_0$ and $j$ identify the
coordinates $t$, $x$, and $y$ on $[-1,1]\times T^2$ with those on $W$.
Hence, we can view $W$ as the further extended time-energy
phase space with the cyclic energy-coordinate $u$. Then $j$ is the graph
of the time-dependent Hamiltonian $K$ in the extended time-energy phase
space $W$. Now it is clear that $v$ is just the Hamiltonian vector field 
of $K$.
\end{Remark}

\begin{Remark}
\labell{rmk:measure-zero}
In the proof of Proposition \ref{prop:K1} we will not require the
Denjoy continuum $\D$ to have zero measure. As a consequence, the
union of characteristics entirely contained in the semi-plug 
can have Hausdorff dimension two because this set is the 
image of $\D$ by a $C^2$-smooth embedding.
\end{Remark}

\section{Proof of Proposition \ref{prop:K1}}
\labell{sec:pfK1}

Recall that $\py+h\px$ is a Denjoy vector field on $T^2$ 
whose choice will be discussed later on and $\D$ is the Denjoy continuum 
for this field. Recall also that $V_0$ is a small, disjoint from $\D$,
neighborhood of $(x_0,y_0)$. Fix a slightly larger neighborhood $V_1$ of
$(x_0,y_0)$ which contains the closure of $V_0$ and is still disjoint
from $\D$. Let $\eps>0$ be sufficiently small.

Proposition \ref{prop:K1} is an immediate consequence of the following

\begin{Proposition}
\labell{prop:K2}
There exists a $C^2$-smooth mapping 
$K\colon [-\eps,\eps]\times T^2\to S^1$ which satisfies
(K1)-(K3) and the requirement

\begin{enumerate}
\item[K4${}^\prime$.]$K=K_0$ for all $t$ and $(x,y)$ in the fixed 
neighborhood $V_1$ of $(x_0,y_0)$.
\end{enumerate}
\end{Proposition}

\begin{proof}[Proof of Proposition \ref{prop:K1}]
Let $K$ be as in Proposition \ref{prop:K2}. We extend this function  
to $[-1,1]\times T^2$ as the linear combination 
$\phi(t)K(t,x,y)+(1-\phi(t))x$, where
$\phi$ is a bump function equal to 1 for $t$ close to 0 and vanishing
for $t$ near $\pm\eps$. Note that this linear combination is well defined
due to (K3). Clearly, the linear combination satisfies (K1)-(K3).
If the range of $t$ for which $\phi(t)=1$ is sufficiently
small it also satisfies (K4).
\end{proof}

\begin{proof}[Proof of Proposition \ref{prop:K2}]~

\emph{Step 1: The extension of $h$ to $[-1,1]\times T^2$.}
Our first goal is to extend $h$ from $T^2$ to 
$H\colon [-1,1]\times T^2 \to \R$ 
smoothly and so that $\px H-\px h$ is of order one in $t$.

\begin{Lemma}
\labell{lemma:ext}
Assume that $\alpha$ is sufficiently close to 1 and 
$h$ is $C^{1+\alpha}$. Then there exists a 
$C^{1}$-function $H\colon [-1,1]\times T^2 \to \R$ such that
\begin{enumerate}
\item[H1.] $H(0,x,y)=h(x,y)$;
\item[H2.] $\px H (t,x,y)= \px h(x,y)+o(t)$ uniformly in $(x,y)$;
\item[H3.] the function $\int_0^t H(\tau,x,y)\,d\tau$ is $C^2$ in $(t,x,y)$.
\end{enumerate}
\end{Lemma}
At this moment only the assertion of Lemma \ref{lemma:ext} is essential and
we defer the proof of the lemma to Section \ref{sec:proofs}.

\begin{Remark}
Since $H$ is only $C^{1}$-smooth, the condition (H2) does
not hold automatically. However, as is easy to see from the proof of the lemma,
one can find an extension $H$ such that $\px H (t,x,y)= \px h(x,y)+o(t^k)$
for any given $k$ and (H3) still holds, provided that $\alpha$ is
sufficiently close 1 (in fact, $k/(k+1)< \alpha <1$).
\end{Remark}

\emph{Step 2: The definition of $K$.} From now on we fix the extension 
$H$, but allow the interval $[-\eps,\eps]$, on which it is considered, 
vary. We will construct the function $K$ of the form
\begin{equation}
\labell{eq:K}
K(t,x,y)=\int_0^t [-H(\tau,x,y) +f(x,y)\tau]\,d\tau+ A(x,y),
\end{equation}
where the ``constant'' of integration $A$ and the correction function
$f$ are chosen so as to make (K1)-(K3) and (K4${}^\prime$) hold. 
Note that
$A$ is actually a function $T^2\to S^1$, whereas $H$ and $f$ are
real valued functions. The main difficulty in the proof below comes
from the combination of the conditions (K1) and (K2).

\emph{Step 3: The auxiliary functions $A$ and $f$.}
Let us now specify the requirements the functions $A$ and $f$ 
have to meet. 

\begin{Lemma}
\labell{lemma:Afg}
There exist a $C^2$-function $A\colon T^2\to S^1$ and 
$C^\infty$-function $f\colon T^2\to \R$ 
satisfying the following conditions:
\begin{enumerate}
\item[A1.] $\px A\geq \eta(\px h)^2$ 
for some constant $\eta>0$ and $\px A$ vanishes exactly on the Denjoy set $\D$;

\item[A2.] there exists an open set $U\subset T^2$, containing $\D$, 
such that $U\cap V_1=\emptyset$ and 
\begin{eqnarray}
\labell{eq:UA}
\px A|_{T^2\ssminus U} &\geq& \const >0,\\ 
\labell{eq:fU}
\px f|_{U} &\geq& 4\eta^{-1}+2;
\end{eqnarray}

\item[A3.] $A$ is a $C^0$-close to $(x,y)\mapsto x$;

\item[A4.] $A(x,y)=x$ for $(x,y)\in V_1$.
\end{enumerate}
\end{Lemma}
This lemma will also be proved in Section \ref{sec:proofs}.

\begin{Remark}
More specifically, the condition (A3) means that for fixed $h$ 
and $V_1$ one can find $A$ arbitrarily $C^0$-close to $(x,y)\mapsto x$
and satisfying other requirements of the lemma.
\end{Remark}

\emph{Step 4: The properties of $K$.}
Let us now prove that $K$ given by \eqref{eq:K}, i.e.,
$$
K=-\int_0^tH\,d\tau +\frac{t^2}{2}f 
+A,
$$
satisfies the requirements of Proposition \ref{prop:K2}, provided
that $\eps>0$ is small enough and $A$ and $f$ are as in Lemma
\ref{lemma:Afg}. The function
$K$ is $C^2$. Indeed, the first term is $C^2$ by (H3). By
Lemma \ref{lemma:Afg}, the next term
is $C^\infty$ and the last term, $A$, is $C^2$.

Condition (K1) is obvious: 
$\pt K|_{t=0} = -H|_{t=0}= -h$ by (H1) and
$\px K|_{\{0\}\times \D} = \px A|_{\D}= 0$ by (A1).

Let us now turn to (K2). We will first show that
\begin{equation}
\labell{eq:inequality}
\px K= -\int_0^t\px H\,d\tau +\frac{t^2}{2}\px f 
+\px A \geq 0
\end{equation}
and then prove that the equality occurs only on $\{0\}\times \D$.

Assume first that $(x,y)\in U$. By (H2) and \eqref{eq:fU}, we have
$$
\px K\geq \px A -t\px h  +2\eta^{-1} t^2
+\left(t^2  + o(t^2)\right).
$$
Obviously,
\begin{equation}
\labell{eq:go4}
t^2  + o(t^2)\geq 0 
\end{equation}
for $(x,y)\in U$ and all $t\in [-\eps, \eps]$, provided that 
$\eps>0$ is small. Hence, to verify \eqref{eq:inequality}, it suffices to 
show that
\begin{equation}
\labell{eq:pfk3}
\px A -t\px h  +2\eta^{-1}t^2  \geq 0.
\end{equation}
By (A1), this follows from
$$
\eta (\px h)^2
-t\px h  +2\eta^{-1}t^2  \geq 0.
$$
Here all the terms are non-negative except, maybe, $-t\px h$.
Hence, it suffices to prove that at least one of the following two 
inequalities holds:
\begin{eqnarray}
\labell{eq:ineq1}
\eta (\px h)^2- t\px h  &\geq& 0,\\
\labell{eq:ineq2}
-t\px h +2\eta^{-1}t^2 &\geq& 0.
\end{eqnarray}
Inequality \eqref{eq:ineq1} holds if (but not only if)
\begin{equation}
\labell{eq:ineq1'}
|t|\leq \eta |\px h|
\end{equation}
and \eqref{eq:ineq2} holds if (but not only if)
\begin{equation}
\labell{eq:ineq2'}
|t|\geq \frac{\eta|\px h|}{2}.
\end{equation}
Clearly, at least one of the 
inequalities \eqref{eq:ineq1'} and \eqref{eq:ineq2'} holds. This proves 
\eqref{eq:inequality} for $(x,y)\in U$.

Assume now that $(x,y)\in T^2\ssminus U$. Then, by (A2) or, more 
specifically, by \eqref{eq:UA}, 
$$
\px K=\px A+O(t)> \const+ O(t)> 0,
$$
when $\eps>0$ is small. Thus (K2) holds for $(x,y)\in (T^2 \ssminus U)$.

To finish the proof of (K2) we need to show that for $(x,y)\in U$ the equality
in \eqref{eq:inequality} implies that $t=0$ and $(x,y)\in 
\D$. Thus, assume that $(x,y)\in U$ and $\px K(t,x,y)=0$. Then
\eqref{eq:go4} and \eqref{eq:pfk3} must become equalities. 
The equality \eqref{eq:go4}
is possible only when $t=0$. Setting $t=0$ in 
the equality \eqref{eq:pfk3}, we conclude that $\px A (x,y)=0$ and hence 
$(x,y)\in \D$ by (A1).

The condition (K3) follows from (A3). Indeed, if $\eps>0$ is small, $K$ is 
$C^0$-close to $A$ which, in turn, is $C^0$-close to $K_0$ by (A3). 

The condition (K4${}^\prime$) need not be satisfied for $K$. By (A3), 
on $[-\eps,\eps]\times T^2$
the function $K$ is $C^0$-close to $K_0$, provided that $\eps>0$ is 
small. Moreover, by (A4), the function $\px K$ is $C^0$-close to 1
and $\py K$ is $C^0$-close to 0 on 
a neighborhood of $[-\eps,\eps]\times \text{closure}(V_1)$, 
for small $\eps>0$. Now it is easy to see that (taking a smaller
$\eps>0$, if necessary) we can modify $K$
near and on $[-\eps,\eps]\times V_1$ so as to keep (K1)-(K3) and make the
new function satisfy (K4${}^\prime$). Indeed, let 
$\phi\colon T^2\to [0,1]$ be a bump function equal to $1$ on $V_1$ and 0
outside of a small neighborhood of $V_1$. Then the linear combination
$x\phi + (1-\phi)K$ still satisfies (K1)-(K3) and also 
(K4${}^\prime$) if $\eps>0$ is small enough.
Note that this linear combination is well defined due to (A3).
\end{proof}

\section{Proofs of Lemmas \ref{lemma:ext} and \ref{lemma:Afg}}
\labell{sec:proofs}

\subsection{Proof of Lemma \ref{lemma:ext}} For $t\neq 0$ and $(x,y)\in T^2$,
set $x_\pm =x \pm t^s/2$ and $y_\pm =y \pm t^s/2$ and define
$$
H(t,x,y)=\frac{1}{t^{2s}}\int_{y_-}^{y_+}\int_{x_-}^{x_+}
h(\xi,\zeta)\,d\xi d\zeta,
$$
where $s$ is an even positive integer to be specified later. Also,
let $H(0,x,y)=h(x,y)$. In other words,
$H(t,x,y)$ is obtained by averaging $h$ over the square with side $t^s$,
centered at $(x,y)$.\footnote{This extension of $h$ by averaging 
is somewhat similar to the one from \cite{kug}.}

\emph{Condition (H2):} 
First note that $H$
is obviously differentiable in $x$ and $y$ for every $t$. Furthermore,
it is easy to see that $H$ satisfies (H2), i.e.,
$\px H=\px h+o(t)$, provided that
\begin{equation}
\labell{eq:s-alpha1}
s\alpha>1.
\end{equation}
Indeed, since $h$ is continuous, we have
$$
\px H(t,x,y) =\frac{1}{t^{2s}}\int_{y_-}^{y_+}
\big( h(x_+,\zeta)-h(x_-,\zeta)\big) \, d\zeta.
$$
By the mean value theorem, 
$h(x_+,\zeta)-h(x_-,\zeta)= t^s \px h(x_0,\zeta)$,
where $x_0$ is some point in $[x_-,x_+]$, depending on $\zeta$.
Since the distance between $(x,y)$ and $(x_0,\zeta)$ does not exceed
$t^s/\sqrt{2}$ and $\px h$ is $\alpha$-H\"older, we have
$|\px h(x_0,\zeta)- \px h(x,y)| \leq \const \cdot \lt t^s \rt^\alpha $ with
$\const$ independent of $(x,y)$. Hence, 
\begin{eqnarray*}
|\px H(t,x,y)-\px h(x,y)|
&\leq &
\frac{1}{t^{s}}\int_{y_-}^{y_+}
\big| \px h(x_0,\zeta)- \px h(x,y)\big| \, d\zeta\\
&\leq& \const \cdot t^{s\alpha},
\end{eqnarray*}
where $\const$ is independent of $(x,y)$. This
proves (H2), provided that \eqref{eq:s-alpha1} holds.

\emph{$C^1$-smoothness of $H$:} We will show that $H$ is $C^1$, provided
that $h$ is $C^1$ and $s>1$. (Note that \eqref{eq:s-alpha1} implies that
$s>1$.) The proof essentially amounts to repeated applications of the 
mean value theorem. However, for the sake of completeness, we give a detailed
argument below.

It is clear that $\px H$ is continuous in all the variables for $t\neq 0$ and
continuous in $x$ and $y$ for all $t$. Its continuity in $t$ at $t=0$
follows from (H2). Moreover, it is easy to see that $\px H\to \px h$ as
$t\to 0$ even if $h$ is just $C^1$. The same reasoning applies to $\py H$.

It remains to show that $\pt H$ exists and is continuous. Again this is
obvious for $t\neq 0$. Using the fact that $h$ is $C^1$,
one can easily check that
$$
|H(t,x,y)-h(x,y)|\leq \const \cdot t^{s}.
$$
This immediately implies that $\pt H\mid_{t=0}=0$ when $s>1$.
Thus, to establish the continuity of $\pt H$ at $t=0$, we need to prove that
$\pt H\to 0$ uniformly in $(x,y)$ as $t\to 0$. A straightforward calculation
shows that
\begin{equation}
\labell{eq:ptH}
\pt H(t,x,y) =
\frac{s}{2t^{s+1}}\int_{y_-}^{y_+}\CI_y(\zeta)\, d\zeta
+
\frac{s}{2t^{s+1}}\int_{x_-}^{x_+}\CI_x(\xi)\, d\xi,
\end{equation}
where 
$$
\CI_y(\zeta)=h(x_+,\zeta)+h(x_-,\zeta)-
\frac{2}{t^s}\int_{x_-}^{x_+} h(\xi,\zeta)\,d\xi
$$
and, similarly,
$$
\CI_x(\xi)=h(\xi,y_+)+h(\xi,y_-)-
\frac{2}{t^s}\int_{y_-}^{y_+} h(\xi,\zeta)\,d\zeta. 
$$
By the mean value theorem, we have
\begin{eqnarray*}
|\CI_y(\zeta)| &=& |h(x_+,\zeta)+h(x_-,\zeta)- 2 h(x_0,\zeta)|\\
               &= & |\px h(x_2,\zeta)(x_+-x_0)-\px h(x_1,\zeta)(x_0-x_-)|\\
               &\leq & \big(|\px h(x_2,\zeta)|+|\px h(x_1,\zeta)|\big)
               \cdot t^s\\
               &\leq & \const \cdot t^s ,
\end{eqnarray*}
where $x_0$, $x_1$, and $x_2$ are some points in $[x_-,x_+]$, depending
$\zeta$ and the constant can be taken independent of $(x,y)$.)
As a consequence, the first term in \eqref{eq:ptH} (whose absolute value
is bounded by $\const\cdot t^{s-1}$) goes to zero as 
$t\to 0$ if $s>1$. A similar argument
shows that the second term also goes to zero. Therefore,
$\pt H\to 0$ as $t\to 0$ uniformly in $(x,y)$, and hence $\pt H$ is
continuous.

\emph{Condition (H3):} 
Let us now prove (H3), i.e., that
$$
F=\int_0^t H\, d\tau
$$
is $C^2$, under some additional constraints on $s$ and $\alpha$. 

It is clear that $F$ is $C^1$. Moreover, the continuity of $\px H$
and $\py H$ implies that $\px \pt F= \px H= \pt \px F$ and
$\py \pt F= \py H= \pt \py F$. Thus these partial derivatives exist
and are continuous.

Furthermore, we claim that the second order partial derivatives of
$F$ in $x$ and $y$ also exist and are continuous, provided that
\begin{equation}
\labell{eq:s-alpha2}
s(1-\alpha)<1.
\end{equation}
Let us examine, for example, $\px\py F$.
Note that $F|_{t=0}=0$, and hence $\px\py F|_{t=0}=0$. Thus we may assume
that $t\neq 0$. Clearly,
\begin{equation}
\labell{eq:pxpyF}
\px \py F
=\px \int_0^t G(\tau,x,y) \,d\tau,
\end{equation}
where
$$
G(\tau, x, y)=
\frac{1}{\tau^{2s}}\int_{x_-}^{x_+}(h (\xi,y_+)-h(\xi,y_-))\,d\xi.
$$
We claim that in \eqref{eq:pxpyF} the integration in $\tau$ and $\px$
can be interchanged. Indeed, 
\begin{eqnarray*}
\bigl|\px G(\tau, x, y)\bigr|
&=&
\frac{1}{\tau^{2s}}
\bigl|
(h (x_+,y_+)-h (x_-,y_+))
 -(h (x_+,y_-)-h (x_-,y_-))
\bigr|\\
&=&
\frac{1}{\tau^{s}}\bigl| \px h (x_2,y_+)-\px h (x_1,y_-)\bigr|.
\end{eqnarray*}
Here $x_1$ and $x_2$ are some points whose distance to $x$ does not
exceed $\tau^s/2$ and the second equality follows from the mean value
theorem. Using the fact that $\px h$ is $\alpha$-H\"older, we obtain
$$
\left|\px G(\tau, x, y)\right|
\leq \frac{\const}{\tau^{s}}(\tau^s)^\alpha
= \frac{\const}{\tau^{s(1-\alpha)}},
$$
where the constant is independent of $x$ and $y$. As a consequence,
if $s(1-\alpha)<1$, i.e., \eqref{eq:s-alpha2} holds,
the integral $\int_0^t \px G(\tau, x, y)\, d\tau$
converges absolutely and uniformly in $(x,y)$. Thus
it follows from \eqref{eq:s-alpha2} that the derivative $\px \py F$ exists
and 
\begin{equation}
\labell{eq:pxpyF2}
\px \py F
=\int_0^t \px G(\tau, x, y)\,d\tau.
\end{equation}
In addition, this implies that
\begin{equation}
\labell{eq:pxpyF-limit}
\px \py F \to 0 \quad\text{as $t\to 0$ uniformly in
$x$ and $y$.}
\end{equation}
Let us prove now that $\px \py F$ is continuous. The above analysis shows
that this derivative is everywhere continuous in $t$ and in $(x,y)$ at $t=0$. 
Hence, we only need to verify its continuity in $(x,y)$ at $t\neq 0$. 
For $t\neq 0$, the integral
\eqref{eq:pxpyF2} can be broken up into two parts: the integral over 
$[0,\delta]$ and the integral over $[\delta, t]$. 
By \eqref{eq:pxpyF-limit}, the first part can be made arbitrarily small
uniformly in $(x,y)$ by choosing $\delta>0$ small. The second part is obviously
continuous in $(x,y)$. This implies that $\px \py F$ is continuous in $(x,y)$.

Other partial derivatives of $F$ in $x$ and $y$ can be dealt with in a 
similar fashion. Hence, to ensure that (H2) and (H3) hold, it
suffices to have $s$ and $\alpha$ satisfy \eqref{eq:s-alpha1} and 
\eqref{eq:s-alpha2} simultaneously. Obviously, for every $\alpha$
sufficiently close to $1$, there exists an even positive integer $s$ 
satisfying these inequalities. This completes the proof of the lemma.

\begin{Remark}
The assumption that $s$ is a positive even integer can be dropped if we
replace $t^s$ by $|t|^s$ in the definition of $H$. Then 
\eqref{eq:s-alpha1} and \eqref{eq:s-alpha2} have a solution $s>0$ if and
only if $1/2<\alpha<1$.
\end{Remark}

\subsection{Proof of Lemma \ref{lemma:Afg}}
\labell{sbsec:lemma-Afg}
\ 

\emph{Step 1: The function $A$.}
The following lemma, which we will prove in Section \ref{sec:pflemma-main},
plays a crucial role in the construction of $A$.

\begin{Lemma}
\labell{lemma:main} There exists a Denjoy vector field $\py +h\px$ which 
is $C^{1+\alpha}$ for all $\alpha \in (0,1)$ and such that
\begin{enumerate}
\item[D1.] $\px h$ vanishes on $\D$;

\item[D2.] $\int_0^x (\px h(\xi,y))^2\,d\xi$ is $C^2$ in $(x,y)$.
\end{enumerate}
\end{Lemma}
 
Assuming Lemma \ref{lemma:main}, let us continue the proof of Lemma 
\ref{lemma:Afg}. The essence of the requirements on $A$ is that
$A$ should be $C^2$, and the derivative $\px A$ should vanish on $\D$ and be 
bounded from below by $\eta(\px h)^2$. If $\eta(\px h)^2$ were not sufficiently
smooth, these conditions would be hard to satisfy. However, since
$\int_0^x \eta (\px h)^2\,d\xi$ is $C^2$ by Lemma \ref{lemma:main}, we may 
simply take
$\eta (\px h)^2$ as $\px A$ with some additional correction terms.
These extra terms are needed to make $A$ into a function $T^2\to S^1$ 
meeting other requirements of Lemma \ref{lemma:Afg}. 

Let us now outline the construction of $A$ omitting some details to be
filled in at the concluding part of the proof (Step 3).
Pick a smooth $C^1$-small non-negative function 
$a\colon T^2\to \R$ which vanishes exactly on $\D$. There exists
a smooth non-negative function $b\colon T^2\to \R$ which vanishes on $\D$
(but not only on $\D$) and such that
\begin{equation}
\labell{eq:b}
(x,y)\mapsto \int_0^x b(\xi,y)\, d\xi
\quad\text{is $C^0$-close to}\quad (x,y)\mapsto x.
\end{equation}
Pick a small $\eta_1>0$ and set\footnote{Throughout the proof we identify 
$T^2$ with $\R^2/\Z^2$.}
$$
A(x,y)=\frac{\int_0^x\left[\eta_1(\px h)^2+a+b\right]\,d\xi}{
\int_0^1\left[\eta_1(\px h)^2+a+b\right]\,d\xi}.
$$
This is a function $T^2\to S^1$. Indeed, $\int_0^1 \px A\, dx=1$
for any $y$ and $A(x,0)=A(x,1)$ for any $x$ by the definition of $A$.
By (D2) and since $a$ and $b$ are smooth, $A$ is $C^2$. By taking
$a$ and $\eta_1>0$ small, we can ensure that $A$ is $C^0$-close to
$(x,y)\mapsto x$, i.e., the requirement (A3) is met.
Also, $A$ obviously satisfies (A1) for some $\eta>0$.

One can construct $b$ in such a way that $b \geq \const$ on the complement 
of some neighborhood $U$ of $\D$, which implies \eqref{eq:UA}, and 
so that $b|_{V_1}=1$. Then on $V_1$, the function $A$ is $C^1$-close to 
$x$, provided
that $\eta_1>0$ is small and $a$ is $C^1$-small. Now it is easy to alter
$A$ on and near $V_1$ so that (A4) is satisfied (i.e., $A(x,y)=x$ on 
$V_1$) and the conditions (A1) and \eqref{eq:UA} still hold.

\emph{Step 2: The function $f$.} First note that it suffices to construct
a function $f$ such that $\px f|_U\geq\const$.  To define such a function $f$,
we pick a smooth function $f'$ such that $f'|_U\geq \const$ and such that 
the mean value of $x\mapsto f'(x,y)$ is zero for every $y$. (This is 
possible if $U$ is sufficiently small.) 
Then $f(x,y)=\int_0^x f'(\xi,y)\,d\xi$ satisfies \eqref{eq:fU}. 

\emph{Step 3: The detailed construction of the neighborhood $U$ and the 
functions $b$ and $f$.} Let us cover $T^2$ by two open overlapping 
cylinders $C_1=S^1\times I_1$ and $C_2=S^1\times I_2$, where $I_1$ and
$I_2$ are two arcs covering the circle $S^1$ with coordinate $y$.

First we describe $b$ and $f$ on $C_1$. For the sake of brevity
let us denote $C_1$ by $C$ and $I_1$ by $I$. Without loss of generality
we may assume that $0\in I$ and $V_1\subset C$. The Denjoy flow gives rise to a
$C^1$-diffeomorphism $\varphi\colon C\to S^1\times I$ which sends 
$\D\cap C$ to a
cylindrical set, i.e., $\varphi(\D\cap C)=\D_0\times I$, where
$\D_0=\D\cap \{y=0\}$. It is easy to see that $\D_0$ can be covered by
a finite collection of disjoint arbitrarily short open intervals 
$\Gamma_1,\ldots, \Gamma_k$. Then $\varphi(\D\cap C)$ is covered by
stripes $\Gamma_i\times I$ and thus $\D\cap C$ is covered by the skewed
stripes $\varphi^{-1}(\Gamma_i\times I)$. 

The intersection of $\varphi^{-1}(\Gamma_i\times I)$ with 
$S^1\times \{y\}$ is an arc whose end-points are $C^1$-functions of $y$. 
For each stripe, let us approximate these functions by $C^\infty$-functions.
If the approximations are accurate enough, the new end-point functions 
still bound non-overlapping skewed open stripes in $C$ which cover $\D\cap C$.
Denote these stripes by $L_1,\ldots, L_k$ and set $U_1=\cup L_i$.

Note that the end-points of $L_i\cap (S^1\times\{y\})$ are smooth functions
of $y$ and that all $L_i$ can be made arbitrarily narrow by taking the
intervals $\Gamma_i$ short. In addition, we can always make $U_1$
disjoint from $V_1$. 

It is not hard to see that there exists a $C^\infty$-function $b$ on $C$ which is 
identically zero on $U_1$ and such that \eqref{eq:b} holds. Indeed, on
$S^1\times \{y\}$ we take a smooth function which is equal to zero on
all arcs $L_i\cap(S^1\times \{y\})$ and has high bumps in between these
arcs. Since the arcs are short, $b$ can be chosen to satisfy
\eqref{eq:b}. This function can be obviously made smooth in $y$ because so
are the end-points of the arcs. In addition, it is easy to see that we can
take $b$ to be equal to $1$ on $V_1$.

The function $f$ is defined in a similar fashion. For example, we can 
take $f'$ equal to $1$ on $U_1$ and, for each $y$, and use the complement of
$U_1\cap (S^1\times \{y\})$ in $S^1\times \{y\}$ to make sure that $f'$ has
zero mean. Then, as we have pointed out above, we set 
$f(x,y)=\int_0^x f'(\xi,y)\,d\xi$.

For the second cylinder $C_2$ the argument is similar. 
The function $b$ 
on $T^2$ is obtained from its counterparts $b_1$ on $C_1$ and $b_2$ 
on $C_2$ by pasting $b_1$ and $b_2$ on $C_1\cap C_2$ using cut-off functions
in $y$. The construction of $f$ is finished in a similar way.
It is clear that there exists
a small neighborhood $U$ of $\D$ (contained in $U_1\cup U_2$) such that
$b|_U=0$ and $\px f|_U=1$. This is the required neighborhood $U$. (Note
that in general we cannot take $U=U_1\cup U_2$.)
The proof of the lemma is completed. 

\section{Proof of Lemma \ref{lemma:main}}
\labell{sec:pflemma-main}
The proof of Lemma ~\ref{lemma:main} is based on the existence of 
a $C^{1+\alpha}$ Denjoy diffeomorphism $\Phi$ such that $(\Phi'-1)^2$ 
is $C^1$.  Therefore, we first outline the construction
of such a Denjoy diffeomorphism, and then proceed with the proof of the lemma.
\subsection{Construction of the Denjoy diffeomorphism}
\labell{sec:constr}
\begin{Lemma} 
\labell{lemma:construction}
There exists 
a Denjoy diffeomorphism $\Phi$ which is $C^{1+\alpha}$ for all 
$\alpha\in (0,1)$ and such that $(\Phi'-1)^2$ is $C^1$.
\end{Lemma}
\begin{proof}[Proof of Lemma ~\ref{lemma:construction}.]
We prove Lemma ~\ref{lemma:construction} in two steps.  First, we define
the required Denjoy diffeomorphism $\Phi$ and show that $\Phi'$ is 
$\alpha$-H\"{o}lder for every $\alpha\in (0,1)$, then
we prove that $(\Phi'-1)^2$ is $C^1$.

\emph{Step 1: Definition of $\Phi$.} In the construction of $\Phi$ we
closely follow the general description of Denjoy maps in 
\cite[Section 12.2]{KH}. Pick $\beta\in (0,1)$ and let
\begin{equation}
\labell{defn:$l_n$}
l_n:=k_\beta(|n|+2)^{-1}(\log\,(|n|+2))^{-1/\beta}
\end{equation}
be the length
of the interval $I_n$ inserted into $S^1$ to ``blow up'' an orbit,  
$a_n$, of an irrational rotation.  Here $k_\beta$ is 
a constant depending on $\beta$ chosen so that 
$\sum_{n\in\integers}\,l_n< 1$.  We emphasize that this choice of $l_n$
is essential in order to make
the series $\sum_{n\in\integers}\,l_n$ converge very slowly which, in turn,
results in a small Denjoy continuum, $S^1 \ssminus 
\bigcup_{n\in \integers}\,Int(I_n)$.  
This slow convergence is the main factor which ensures that 
$(\Phi'-1)^2$ is $C^1$ and the second assertion (D\ref{itm:d2}) 
of Lemma ~\ref{lemma:main} holds, as it will become clear later on.     

To construct a Denjoy diffeomorphism $\Phi$, it suffices to define the 
derivative $\Phi'$, since $\Phi$ is then obtained by integration.    
Let $\varphi\colon[0,1]\rightarrow\reals$ be a bump function 
satisfying $\int_0^1 \varphi(x)\,dx=1$.  Define 
the smooth function 
$$
\varphi_n(x):=c_n\,\varphi\big((x-a_n)/{l_n}\big)
$$ 
on the interval $I_n=[a_n, a_n+l_n]$,
where $c_n=(l_n-l_{n+1})/{l_n}$, and note that
$\int_{I_n}\varphi_n(x)\,dx=c_n\,l_n$.  Finally, let  
\[\Phi'(x) = \begin{cases}
               1& \text{for $x\,\in S^1 \ssminus 
                            \bigcup_{n\in \integers}\,I_n$},\\
  1+\varphi_n(x)& \text{for $x\,\in \,I_n$}.
              \end{cases}\]
This completes the construction of $\Phi$. It is well known, \cite{KH}, 
and easy to see that $\Phi$ is $C^{1+\alpha}$ for any $\alpha\in (0,1)$.
(Moreover, one can show that
$|\Phi'(x)-\Phi'(x_0)|\leq \const |x-x_0|\big|\log |x-x_0|\big|^{1/\beta}$ 
for any $x$ and $x_0$ in $S^1$.) For what follows, we only need that
$\Phi$ is $C^{1+\alpha}$ for some $\alpha\in (1/2,1)$ and also some estimates
on the $C^1$-norm of $(\Phi'-1)|_{I_n}$ which result from \eqref{defn:$l_n$}.

Let us now list some properties of $c_n$ and $l_n$ that we will use later on:  

First, we note that, as is true for any Denjoy map,
$$
c_n\to 0\quad\text{as}\quad n\to\infty.
$$
In fact, $c_n=O(1/n)$.

Furthermore, \eqref{defn:$l_n$} guarantees\footnote{This is the 
the main point in the proof where the specific choice of $l_n$ made 
above is essential.}
that
\begin{equation}
\labell{est:$c_n^2/l_n$} 
 \fr{c_n^2}{l_n}
\to 0\quad\text{as}\quad n\to \infty.
\end{equation}
Indeed,
$$
\fr{|c_n|}{l_n}=
\fr{ \left|l_n-l_{n+1}\right|}{l_n^2}=
O \bigg( \big( \log\,\lt |n|+2 \rt \big)^{1/\beta}\bigg), 
$$
as can be seen by expanding the left hand side in 
$|n|^a\,\big( \log\,\lt |n|+2 \rt \big)^b$ for $a\leq 0$, and $b\geq 0$. 
Thus, 
$$
 \fr{c_n^2}{l_n}= 
 l_n\,\fr{c_n^2}{l_n^2}=
 \fr{O\big( \big( \log\,(|n|+2) \big) ^ {1/ \beta} \big)}
 {\lt |n|+2 \rt }
\to 0\quad\text{as}\quad n\to \infty.
$$

We finish this discussion by
establishing the following estimates which will be used in the rest
of the proof\footnote{Throughout the rest of the proof $\|~\|$ denotes
the $\sup$-norm on $I_n$.}:
\begin{eqnarray}
\labell{ineq:Phi}
\left\|\left(\Phi'-1\right)|_{I_n}\right\| &=& O(|c_n|)\to 0\\
\labell{ineq:der.Phi}
  \big\|\px\left(\Phi'-1\right)|_{I_n}\big\| &=& O(|c_n|/l_n)\\
\labell{ineq:der.square}  
  \big\|\px\left(\Phi'-1\right)^2|_{I_n}\big\| &=& O(c_n^2/l_n)\to 0
\end{eqnarray}
To prove these estimates, we first recall that $(\Phi'-1)|_{I_n}=\varphi_n$.
Then, since $\|\varphi_n\|=|c_n|\cdot \|\varphi\|$, we have 
\eqref{ineq:Phi}. The second estimate, \eqref{ineq:der.Phi}, is proved as 
follows:
$$
  \left\|\px\left(\Phi'-1\right)|_{I_n}\right\|=
  \left\|\px\varphi_n\right\| \leq
  \|\varphi'\|\fr{|c_n|}{l_n}. 
$$
Finally, \eqref{ineq:der.square} is a consequence of the first two estimates
and \eqref{est:$c_n^2/l_n$}.

\emph {Step 2: Proof that $(\Phi'-1)^2$ is $C^1$.} Let 
$\D_0:=S^1 \ssminus \bigcup_{n\in \integers}\,Int(I_n)$
denote the Denjoy continuum.  
Observe that, since on each $I_n$ the function $\varphi_n$ is 
smooth, $\Phi'-1$ as well as $(\Phi'-1)^2$ are also smooth on $I_n$.  
Hence, we  need to prove that for $x_0\in \D_0$, $(\Phi'-1)^2$ is 
differentiable,  
its derivative at $x_0$ is zero, and $\px\left(\Phi'-1\right)^2(x) 
\to 0$ as $x \to x_0$. 

Recall that $\Phi'-1$ is $\alpha$-H\"{o}lder continuous with
$\alpha > 1/2$ and $(\Phi'-1)^2\equiv0$ on $\D_0$. It readily
follows that $(\Phi'-1)^2$ is differentiable, and its derivative is zero 
on $\D_0$.

To finish the proof,
it remains to show that $\px\left(\Phi'-1\right)^2(x)\to 0$
as $x\to x_0\in \D_0$. Let $x_k$ be a sequence in $S^1 \ssminus 
\D_0$ converging to $x_0$.  
Since $\D_0$ is nowhere dense, there exists a sequence of intervals, 
$I_{n_k}$, such that $x_k \in I_{n_k}$ for $k\in \mathbb N$. 
Then, by \eqref{ineq:der.square}, 
$\left|\px\left(\Phi'-1\right)^2\,(x_k)\right|
\rightarrow 0$ as $k \rightarrow \infty$, and this, together with  
$\px\left(\Phi'-1\right)^2(x_0)=0$, proves the assertion.
\end{proof}

\begin{Remark}
\labell{rmk:smooth}
The Denjoy map defined by \eqref{defn:$l_n$} is essentially as smooth
as a Denjoy map can be made, up to using functions growing slower than 
logarithms, e.g., iterations of logarithms. The next significant 
improvement in smoothness would be to have $\log \Phi'$ of bounded variation
or satisfying the Zygmund condition which is impossible; see \cite{HS,JS,KH}.
\end{Remark}

Now we are in a position to prove Lemma ~\ref{lemma:main}
which asserts:
\emph{There exists a Denjoy vector field $\py +h\px$ which 
is $C^{1+\alpha}$ for all $\alpha \in (0,1)$ and such that}
 \begin{list}{\rm D\arabic{cond}.}{\usecounter{cond}}
   \item \label{itm:d1} \emph{$\px h$ vanishes on $\D$;}
   \item \label{itm:d2} \emph{$\int_0^x\,(\px h(\xi,y))^2\,d\xi$ is
$C^2$ in $(x,y)$.}
 \end{list}
To prove this lemma we show that the Denjoy vector field $\py+h\px$
on $T^2$ for $\Phi$ described above satisfies (D\ref{itm:d1}) and 
(D\ref{itm:d2}). The proof of (D\ref{itm:d1}) is straightforward and
based on the explicit formula for $h$. 
The proof of (D\ref{itm:d2}) is divided into two parts. In the first part
(Section \ref{sbsbsec:$h'_x$is$C^1$}), we show that  
$\lt \px h \rt ^2 $ is $C^1$, which obviously means that
$\px \int_0^x\,(\px h(\xi,y))^2\,d\xi$ is $C^1$.
In the second part (Section \ref{subsec:py}), we show that
$\py \int_0^x\,(\px h(\xi,y))^2\,d\xi$ is $C^1$.  These two results
imply (D\ref{itm:d2}).

\subsection{Explicit formula for $h$ and the proof of (D\ref{itm:d1}).}
\labell{sec:h}
First we give explicit formulas for $h$ and $\px h$, and 
fix some notations. Let
\begin{equation}
\labell{defn:circle diff}
\Phi_y(x)=\lt 1-\delta \lt y \rt \rt  x+
         \delta \lt y \rt  \Phi \lt x \rt, 
\end{equation}
where $\delta\colon[0,1]\rightarrow[0,1]$ is a non-negative, increasing, 
smooth function which is $0$ for $y$ close to $0$ and $1$ for 
$y$ close to $1$. Then the $x$-component $h$ of a Denjoy vector field
can be expressed as 
$$
h(x,y)=\lt\py\Phi_y\circ\Pin\rt \lt x \rt,
$$ 
where the function $\Pin\lt x \rt$ 
is the inverse of $\Phi_y(x)$ in the $x$-variable.
 
Analyzing the smoothness of these functions, we first observe that
$\Phi_y(x)$ is clearly $C^{1+\alpha}$ in $(x,y)$. 
Moreover, $\Phi_y(x)$ is $C^\infty$ for $x\not\in \D_0$. Furthermore,
$\Pin(x)$ is also $C^{1+\alpha}$ in $(x,y)$. To see this note that
by the implicit function theorem $\Pin(x)$ is $C^1$ in $(x,y)$ and 
$$
\py \Pin(x)=-\frac{\py \Phi_y\lt \Pin(x)\rt}{\px \Phi_y\lt \Pin(x)\rt},
$$
where the denominator is bounded away from zero. Now it readily follows
that $\py \Pin (x)$ is $C^\alpha$ in $(x,y)$, for the numerator is 
$C^{1+\alpha}$ and the denominator is $C^\alpha$. A similar argument shows
that $\px \Pin (x)$ is $C^\alpha$ in $(x,y)$.

As a consequence, $h$ is $C^{1+\alpha}$,  and hence
$$\px h \lt x,y \rt = \delta'(y)
       \lt \Phi' \circ \Pin  \lt x \rt -1 \rt
        \px \Pin  \lt x \rt $$
is $C^\alpha$.

Finally, for a fixed $y\in[0,1]$, keeping the notation from
Section \ref{sbsec:lemma-Afg},   
let $\D_y := \D \cap \{ y \} $. Thus, $\D_y = \Phi_y \lt \D_0 \rt$.   

\noindent{\em Proof of (D\ref{itm:d1}).} Let $(x,y)\in \D$, i.e., 
$x \in \D_{y} = \Phi_y \lt \D_0 \rt $.   
Thus, $\Pin \lt x \rt \in \D_0$. Since $\Phi'-1\equiv0$ on 
$\D_0$, we conclude that $\px h(x,y)=0$, i.e., $\px h$ vanishes on $\D$.

\subsection {Proof of (D\ref{itm:d2}), Part I: $\lt\px h\rt^2$ is $C^1$.}
\label{sbsbsec:$h'_x$is$C^1$}
Note that the existence and continuity of the partial derivatives of 
$\lt \px h \rt ^2 $ is non-trivial only at the points of $\D$.

First, observe that both of the partial  derivatives $\px \lt\px h\rt^2$ 
and $\py \lt\px h\rt^2$ exist and vanish at $(x,y)\in \D$. This follows 
immediately from the facts that $\px h$ is
$\alpha$-H\"{o}lder continuous with $\alpha>1/2$
and $\px h$ vanishes on $\D$ by (D1). 

To examine the continuity of $\px \lt \px h \rt ^2 $ and 
$\py \lt \px h \rt ^2 $ we adopt a new notation 
for $\lt \px h \rt ^ 2$. 
Fix $y\in [0,1]$ and let $F_y\colon S^1 \rightarrow \reals $  be the function 
defined by 
\begin{equation}
\labell{eq:$F_y$}
F_y(\xi) = \big( \Phi' \lt \xi \rt -1 \big) ^2
       \fr{\big( \delta'\lt y \rt  \big) ^2}
      {\big( 1 + \delta \lt y \rt \lt \Phi' \lt \xi \rt -1 \rt \big) ^ 2}. 
\end{equation}
Then 
\begin{equation}
\labell{eq:h-F}
\lt \px h(x,y) \rt ^2 = F_y \circ \Pin  \lt x \rt.
\end{equation}

It follows that $F_y$ vanishes on $\D_0$ for every $y$.
The function $F_y$ is clearly differentiable since 
$F_y \lt \xi \rt = \lt \px h \rt ^2 \lt \Phi_y \lt \xi \rt,y \rt$, where
$\lt \px h \rt ^2$ is differentiable as is shown above and 
$\Phi_y \lt \xi \rt$ is $C^{1+\alpha}$ as proved in Section
\ref{sec:h}.  Furthermore, $\pxi F_y$ and 
$\py F_y$ are both zero on $\D_0$, for the partial derivatives of 
$\lt \px h \rt ^2 $ vanish on $\D$.

To prove that $\lt \px h \rt ^2$ is $C^1$ in $(x,y)$, it suffices to show 
that $F_y(\xi)$ is $C^1$ in $(\xi,y)$. (Indeed, $\Pin$ is $C^{1+\alpha}$ and 
\eqref{eq:h-F} implies that $\lt \px h \rt^2$ is $C^1$ if $F_y$ is $C^1$.)
Thus, it remains to prove that $\pxi F_y $ and $\py F_y$ are continuous.  

\emph{Continuity of $\py F_y(\xi)$.} This follows immediately from
\eqref{eq:$F_y$} since $\delta$ is $C^\infty$-smooth.

\emph{Continuity of $\pxi F_y(\xi)$.} First note that a straightforward 
calculation using \eqref{eq:$F_y$} shows that 
$$\pxi F_y(\xi) = \big( \delta' \lt y \rt \big) ^2
                \fr{ \big( 1+\delta \lt y \rt \lt \Phi'-1 \rt \big)
                     \,\pxi \lt \Phi'-1 \rt ^2 - 2 \,\delta \lt y \rt 
                     \lt \Phi'-1 \rt ^2 \pxi \lt \Phi'-1 \rt }
                   {\big( 1+\delta \lt y \rt \lt \Phi'-1 \rt \big) ^3}$$ 
on $S^1 \ssminus \D_0=\bigcup_{n\in \integers}\,Int(I_n)$, and, 
as discussed above, $\pxi F_y(\xi)=0$ on $\D_0$. It follows immediately
that $\pxi F_y(\xi)$ is continuous in $y$ for every $\xi$.
 
Clearly,  $\pxi F_y(\xi)$ is continuous in $\xi$ on the complement of
$\D_0$ for every fixed $y$. Let us show the continuity at 
$(\xi,y)$ with $\xi \in \D_0$. Note that
the denominator in the expression for $\pxi F_y$ is bounded away from
zero. Using the estimates \eqref{ineq:Phi}, \eqref{ineq:der.Phi}, and 
\eqref{ineq:der.square},  it is easy to
see that the asymptotic behavior of $\|\pxi F_y|_{I_n}\|$ as $n\to\infty$
is determined by $\|\pxi \lt \Phi'-1 \rt ^2\|$, i.e.,
\begin{equation}
\labell{eq:added}
\big\| \pxi F_y|_{I_n} \big\|= O({c_n^2}/l_n)\to 0.
\end{equation}
Arguing as in the proof of the fact that
$\lt \Phi'-1 \rt ^2$ is continuously differentiable (see Section 
\ref{sec:constr}), we conclude that $\pxi F_y(\xi)$ is continuous.

This finishes the proof that $F_y$, and hence $ \lt \px h \rt ^2$,
is $C^1$. 

\subsection {Proof of (D\ref{itm:d2}), Part II: 
$\py \int_0^x\,(\px h(\xi,y))^2\,d\xi$ is $C^1$.} 
\labell{subsec:py}
First let us write the function $\int_0^x\,(\px h(\xi,y))^2\,d\xi$ in a form
more convenient for our analysis. Setting $\xi=\Phi_y(\eta)$, we obtain
\begin{eqnarray*}
    \int_0^x\,(\px h(\xi,y))^2\,d\xi
&=& \int_0^x\, F_y \lt \Pin \lt \xi \rt \rt \,d\xi\\
&=& \int_{\Pin(0)}^{\Pin(x)}\, F_y(\eta)\, \peta \Phi_y(\eta)\,d\eta.
\end{eqnarray*}
Hence, define
\begin{equation}
\labell{eq:$G_y$}
G_y(u):= \int_0^u\, F_y(\eta) \,\peta \Phi_y(\eta)\,d\eta.
\end{equation}
Then
\begin{eqnarray*}
    \int_0^x\,(\px h(\xi,y))^2\,d\xi
&=& \int_{\Pin(0)}^{\Pin(x)}\, F_y(\eta)\, \peta \Phi_y(\eta)\,d\eta\\
&=& \int_0^{\Pin(x)}\, F_y(\eta)\, \peta \Phi_y(\eta)\,d\eta
    -\int_0^{\Pin(0)}\, F_y(\eta)\, \peta \Phi_y(\eta)\,d\eta\\
&=& G_y\circ\Pin(x) - G_y\circ\Pin(0)
\end{eqnarray*}

Thus, our goal is to prove that $\py [G_y\circ\Pin(x)]$ is a $C^1$-function.
We do this in two steps: first, we show that the function 
$G_y(u)$ is $C^2$ and then, using this result, we prove that 
$\py [G_y\circ\Pin(x)]$ is a $C^1$-function.

\subsubsection{Step 1: Proof that $G_y(u)$ is $C^2$.}
\labell{sec:Gy}
Let us show that both of the partial derivatives $\pu G_y(u)$ and
$\py G_y(u)$ are $C^1$.

\noindent{\em Proof that $\pu G_y(u)$ is $C^1$.}
Let $\tF_y(u):=\pu G_y(u)$. Thus, by \eqref{eq:$G_y$}, 
\begin{equation}
\labell{eq:tF}
\tF_y(u) = F_y(u) \,\pu \Phi_y(u).
\end{equation}

First, let us consider $\py\pu G_y=\py \tF_y$. 
For $u\in S^1\ssminus\D_0$, the function $\tF_y(u)$ is smooth.  Hence, 
as long as $u\in S^1\ssminus\D_0$, the derivative $\py \tF_y$ exists (and is
continuous). For $u_0$ in $\D_0$, the partial derivative 
$\py\tilde{F_y}(u_0)$ exists and is zero.  The reason is that 
$\tilde{F_y}(u_0)=0$ for all $y\in[0,1]$.  Furthermore, 
$\py \tF_y(u)$ is continuous in $u$ and smooth in $y$,
i.e., $\py \tF_y(u)$ is infinitely differentiable in $y$ and every derivative
is continuous in $(u,y)$ as immediately follows from
\eqref{defn:circle diff} and \eqref{eq:$F_y$}. This proves the continuity
of $\py\pu G_y$.

Let us now focus on the partial derivative $\pu^2 G_y=\pu\tF_y$. 
As before, this partial derivative obviously exists when 
$u\in S^1\ssminus\D_0$. Furthermore, we claim that 
$\pu \tilde{F_y}(u_0)$ exists and is zero for any $u_0\in\D_0$. 
To see this, recall that as we proved in Section
\ref{sbsbsec:$h'_x$is$C^1$}, $F_y(u_0)=0$ and $\pu F_y(u_0)=0$ for all 
$u_0\in\D_0$. Hence, 
\begin{eqnarray*}
   \pu \tilde{F_y}(u_0)
&=&\pu \big( F_y(u)\, \pu \Phi_y(u) \big)|_{u=u_0}\\
&=&\lm_{u \rightarrow u_0}{\fr{F_y(u)\, \pu \Phi_y(u)-
                           \overbrace{F_y(u_0)}^0 \,\pu \Phi_y(u_0)}
                               {u-u_0} }\\
&=&\lm_{u \rightarrow u_0}{ \fr{F_y(u)-F_y(u_0)}{u-u_0}\,\pu \Phi_y(u) }\\
&=&\underbrace{\pu F_y(u_0)}_0\,\pu \Phi_y(u_0)\\
&=& 0.
\end{eqnarray*}
                               
To show that $\pu \tF_y$ is continuous, we first express $\pu \tF_y$ on 
each $I_n$ as follows
$$\pu \tF_y(u) =\underbrace{\pu F_y(u)}_{O({c_n^2}/l_n)} \, 
                \underbrace{\pu \Phi_y (u)}_{O(1)} 
               +\underbrace{F_y(u)}_{O(c_n^2)}\,
                \underbrace{\pu^2 \Phi_y (u)}_{O(|c_n|/l_n)},$$
where the braces indicate asymptotic behavior as $|n|\to\infty$.
The estimate $\|\pu \tF\|=O(c_n^2/l_n)$ has been established above,
see \eqref{eq:added}; the estimate $\|\pu \Phi_y\|=O(1)$ follows from
the definition of $\Phi$ (see \eqref{defn:circle diff}) and \eqref{ineq:Phi};
the estimate $\| F_y\|=O(c_n^2)$ is a consequence of the definition of
$F_y$ (i.e., \eqref{eq:$F_y$}) and \eqref{ineq:Phi}. Finally,
$\|\pu^2 \Phi_y\|=O(|c_n|/l_n)$ results from \eqref{defn:circle diff}
and \eqref{ineq:der.Phi}.

Now it is clear that 
$\|\pu \tF_y(u)_{I_n}\| = O({c_n^2}/l_n)$.  Since, by 
(\ref{est:$c_n^2/l_n$}), ${c_n^2}/l_n\rightarrow 0$ 
as $|n|\rightarrow \infty$, $\pu\tF_y(u)$ can be shown to be 
continuous in a fashion similar to the cases discussed before. 

This completes the proof of the fact that $\pu G_y(u)=\tF_y(u)$ is $C^1$.

\noindent{\em Proof that $\py G_y(u)$ is $C^1$.}
Note that, since $\py \tF_y(\eta)$ is continuous in $(\eta,y)$ and its domain 
is compact, the functions $\py \tF_y$ converge uniformly to 
$\py \tF_y |_{y=y_{0}}$ as $y \rightarrow y_0$ for any $y_{0}\in[0,1]$.  
Thus, we have
$$\py G_y(u)=\py \int_0^u\, \tF_y(\eta) \,d\eta
                     =\int_0^u\, \py \tF_y(\eta) \,d\eta.$$
This implies that $\pu\py G_y(u)$ is continuous, for
$\pu\py G_y(u)=\tF_y(u)$ is continuous (in fact, $C^1$).
To show that $\py^2 G_y$ is continuous we recall that $\tF_y(u)$ is 
infinitely differentiable in $y$ and every derivative is continuous in
$(u,y)$. Hence, as above, the integration and differentiation can be 
interchanged, and
$$
\py^2 G_y (u)=\int_0^u\, \py^2 \tF_y(\eta) \,d\eta
$$
is continuous because the integrand is contiunous. This completes Step 1.

\subsubsection{Step 2: Proof that $\py [ G_y\circ\Pin(x) ]$ is $C^1$.}
We first write this partial derivative explicitly as follows
$$\py \big( G_y\circ\Pin(x) \big)=\pu G_y \big( \Pin(x) \big)\,\py \Pin(x) 
                                  +\py G_y \big( \Pin(x) \big).$$
The second term of the sum on the right hand side is $C^1$ because 
it is the composition of $C^1$ and $C^{1+\alpha}$ functions. Thus, 
we focus on the first summand which is 
\begin{eqnarray*}
\labell{eq:final}
\pu G_y \big( \Pin(x) \big)\,
\py \Pin(x) &=& \tilde{F_{y}} \big( \Pin(x) \big)\,\py \Pin(x)\\
&=& F_{y} \big( \Pin(x) \big)\,\pu \Phi \big( \Pin(x) \big)\,\py \Pin(x),
\end{eqnarray*}
where the last equality follows from \eqref{eq:tF}.
The product of the last two terms can be further simplified. Applying $\py$ to
the identity $\Phi_y\big(\Pin(x)\big)\equiv x$, we obtain
$$
\pu \Phi_y \big( \Pin(x) \big)\,\py \Pin(x)
+ \lt \py \Phi_y \rt \big( \Pin(x) \big) =0,
$$ 
and hence
$$
\tilde{F_{y}} \big( \Pin(x) \big)\,\py \Pin(x)
=-[F_{y} \,\py \Phi_y]\circ  \Pin(x) .
$$
Recall that $F_y(u)$, $\py \Phi_y(u)$ and $\Pin(x)$ are all $C^1$-functions.
Hence, the left hand side is also $C^1$.

This concludes Step 2 and hence the proof of the fact that 
$\py \int_0^x\,(\px h(\xi,y))^2\,d\xi$ is $C^1$.

\begin{Remark}
Note that the norms of $(\Phi'-1)^2|_{I_n}$ and 
$\pu F_y(u)|_{I_n}$ 
and $\pu \tF_y(u)|_{I_n}$ converge to zero only as $O(c_n^2/l_n)$. 
(One can also show that the same is true for 
the $\px$- and $\py$-partial derivatives of $G_y\circ \Phi^\inv_y(x)$.) A 
faster rate of convergence, e.g., $O(|c_n|^3/l_n)$, would be likely to result
in an ``unacceptably'' smooth Denjoy map and vector field.
\end{Remark}


\begin{thebibliography}{ABCD}

\bibitem[Gi1]{gi:seifert}
V. L. Ginzburg, An embedding $S^{2n-1}\to\reals^{2n}$, $2n-1\geq 7$,
whose Hamiltonian flow has no periodic trajectories, 
{\em IMRN}, (1995), no. 2, 83--98.

\bibitem[Gi2]{gi:seifert97}
V. L. Ginzburg, A smooth counterexample to the Hamiltonian Seifert
conjecture in $\reals^6$, {\em IMRN}, (1997), no. 13, 641--650.

\bibitem[Gi3]{gi:bayarea}
V. L. Ginzburg, Hamiltonian dynamical systems without periodic orbits, in
\emph{Northern California Symplectic Geometry Seminar}, 35--48, 
Amer. Math. Soc. Transl. Ser. 2, vol. 196, Amer. Math. Soc., 
Providence, RI, 1999.

\bibitem[Gi4]{gi:barcelona}
V. L. Ginzburg, The Hamiltonian Seifert conjecture: examples and open
problems, math.DG/0004020; to appear in the \emph{Proceedings of the 
Third ECM, Barcelona, 2000}, Birkh\"auser.

\bibitem[GG]{GG}
V. L. Ginzburg and B. Z. G\"urel, On the construction of a 
$C^2$-counterexample to the Hamiltonian Seifert Conjecture in $\R^4$,
Preprint 2001.

\bibitem[He1]{herman-fax}
M.-R. Herman, Fax to Eliashberg, 1994.

\bibitem[He2]{herman}
M.-R. Herman, Examples of compact hypersurfaces in $\R^{2p}$, $2p\geq 6$, 
with no periodic orbits, in \emph{Hamiltonian systems with three or more
degrees of freedom}, C. Simo (Editor), NATO Adv. Sci. Inst. Ser. C, Math.
Phys. Sci., vol. 533, Kluwer Acad. Publ., Dordrecht, 1999.

\bibitem[HZ1]{ho-ze:per-sol}
H. Hofer and E. Zehnder, Periodic solution on hypersurfaces and a
result by C. Viterbo, {\em Invent. Math.}, {\bf 90} (1987), 1--9.

\bibitem[HZ2]{ho-ze:book}
H. Hofer and E. Zehnder, {\em Symplectic invariants and Hamiltonian
dynamics}, Birkh\"{a}user, Advanced Texts; Basel-Boston-Berlin, 1994.


\bibitem[HS]{HS}
J. Hu and D. Sullivan, Topological conjugacy of circle diffeomorphisms, 
\emph{Ergodic Theory Dynam. Systems}, \textbf{17} (1997), 173--186. 

\bibitem[JS]{JS}
M. Jakobson and G. \'Swi\c atek, One-dimensional maps, Preprint.

\bibitem[KH]{KH}
A. Katok and B. Hasselblatt, 
\emph{Introduction to the modern theory of dynamical systems}, Encyclopedia of 
Mathematics and its Applications, 54. 
Cambridge University Press, Cambridge, 1995. 

\bibitem[Ke]{ely:example}
E. Kerman, New smooth counterexamples to the Hamiltonian Seifert conjecture,
Preprint 2001, math.DG/0101185; to appear in the \emph{Journal of Symplectic
Geometry}.

\bibitem[KuG]{kug}
G. Kuperberg, A volume--preserving counterexample to the Seifert 
conjecture, {\em Comment. Math. Helv.}, {\bf 71} (1996), 70--97.

\bibitem[KuGK]{kugk}
G. Kuperberg and K. Kuperberg, Generalized counterexamples to the 
Seifert conjecture, {\em Ann. Math.}, {\bf 144} (1996), 239--268.

\bibitem[KuK1]{kuk}
K. Kuperberg, A smooth counterexample to the Seifert conjecture in
dimension three, \emph{Ann. Math.}, \textbf{140} (1994), 723--732.

\bibitem[KuK2]{kuk:icm}
K. Kuperberg,  Counterexamples to the Seifert conjecture,
\emph{Proceedings of the International Congress of Mathematicians},
Vol. II (Berlin, 1998). Doc. Math. (1998) Extra
Vol. II, 831--840.


\bibitem[KuK3]{kuk:notices}
K. Kuperberg, Aperiodic dynamical systems, \emph{Notices Amer. Math.
Soc.}, \textbf{46} (1999), 1035--1040.

\bibitem[McDS]{mcduff-sal}
D. McDuff and D. Salamon, \emph{Introduction to symplectic topology}, Oxford
Mathematical Monographs. Oxford University Press, New York, 1995. 

\bibitem[Sc]{schweitzer}
P. A. Schweitzer, Counterexamples to the Seifert conjecture and
opening closed leaves of foliations, \emph{Ann. Math.},
{\bf 100} (1970), 229--234. 

\bibitem[St]{str}
M. Struwe, Existence of periodic solutions of Hamiltonian
systems on almost every energy surfaces, {\em Bol. Soc. Bras. Mat.},
{\bf 20} (1990), 49--58.

\end{thebibliography}
\end{document}